\documentclass[10pt,twoside]{article}

\date{}
\usepackage{amsmath}
\usepackage{amssymb}
\usepackage{amsthm}

\newcommand{\C}{\mathbb{C}}

\newcommand{\g}{\mathfrak{g}}

\newcommand{\lk}{\mathfrak{k}}

\newcommand{\m}{\mathfrak{m}}

\newcommand{\bdm}{\begin{displaymath}}
\newcommand{\edm}{\end{displaymath}}

\theoremstyle{definition}

\newtheorem{lem}{Lemma}
\newtheorem{thm}{Theorem}
\newtheorem{defn}{Definition}

\newtheorem{rem}{Remark}

\newtheorem{axiom}{Corollary}

\title{Adding a uniton via the DPW method}
\author{N. Correia and R. Pacheco}

\begin{document}

\maketitle

\emph{Departamento de Matem\'{a}tica, Universidade da Beira
Interior, Rua Marqu\^{e}s d'\'{A}vila e Bolama 6201-001
Covilh\~{a} - Portugal}

\emph{email: ncorreia@mat.ubi.pt; rpacheco@mat.ubi.pt }

\begin{abstract}
In this paper we  describe how the operation of adding a uniton
arises via the DPW method of obtaining harmonic maps into compact
Riemannian symmetric spaces out of certain holomorphic one forms.
We  exploit this point of view to investigate which unitons
preserve finite type property of harmonic maps. In particular, we
 prove that the Gauss bundle of a harmonic map of finite type
into a Grassmannian is also of finite type.
\end{abstract}



\section{Introduction}

In \cite{Uh}, Uhlenbeck observed that, in the setting of harmonic
maps of $\mathbb{R}^2$ into a compact Riemannian symmetric space
$G/K$, the harmonic map equations amount to the flatness of a
family of connections depending on an auxiliary parameter
$\lambda\in S^1$. This zero-curvature formulation yields an action
of a certain loop group on the space of harmonic maps; underlying
this \emph{dressing} action is the existence of Iwasawa-type
decompositions of the loop groups and loop algebras concerned (cf.
\cite{bG,BP2,DPW,GO,Uh}.)

In \cite{DPW}, the authors introduced a systematic procedure of
obtaining harmonic maps into a compact Riemannian symmetric space
$G/K$ from holomorphic $1$-forms with values in the subspace
\begin{equation}\label{prt}
 \Lambda_{-1,\infty}=\left\{\xi\in\Lambda
\mathfrak{g}^\C\, |\,\,\mbox{$\lambda \xi$ extends holomorphically
to $ |\lambda|<1$}\right\}, \end{equation} where $\g^\C$ denotes
the Lie algebra of $G^\C$, and proved that this correspondence
between $\Lambda_{-1,\infty}$-valued holomorphic $1$-forms $\mu$,
the \emph{holomorphic potentials}, and harmonic maps is
equivariant with respect to loop group actions by dressing and by
gauge transformations on the $1$-forms $\mu$.

Another well known operation for generating new harmonic maps from
a given one was introduced by Uhlenbeck  \cite{Uh} and is called
{\it{adding a uniton}}. In the case of harmonic maps into
Grassmannians, this procedure corresponds to the \emph{forward
replacements} and \emph{backward replacements} of Burstall and
Wood \cite{BW}.

In this paper we describe how the operation of adding a uniton
arises via gauge transformations on the holomorphic $1$-forms
$\mu$ and we exploit this point of view to investigate which
unitons preserve \emph{finite type} property. Recall that harmonic
maps of finite type are obtained  by integrating a pair of
commuting Hamiltonian vector fields on certain finite-dimensional
subspaces of loop algebras (cf. \cite{BFP,BP}) and play a
fundamental role in the genus one case, i.e., harmonic maps from a
two-torus -- for example, it was shown in \cite{BFP} that any
non-conformal harmonic map of a two-torus into a rank one
symmetric space $G/K$ is of finite type. In particular, we  prove
that the Gauss bundle of a harmonic map of finite type into a
Grassmannian is also of finite type. Finally we show that these
unitons preserving finite type can be added by taking a limit of
dressing transformations as in the \emph{completion} procedure
studied by Bergvelt and Guest in \cite{bG}. To obtain these
results one has to enlarge the loop groups involved in DPW
procedure, as in \cite{BP2}, and the space of holomorphic
potentials $\mu$ (see Section \ref{hol}).


\section{Extended solutions} We start by summarizing briefly the
relevant definitions and results concerning the well known
correspondence between harmonic maps and \emph{extended
solutions}, referring the reader to the seminal paper \cite{Uh}
for details.

\vspace{.15in}

Let $G$ be a compact (connected) semisimple  matrix Lie group,
with identity $e$ and Lie algebra $\mathfrak{g}$. Equip $G$ with a
bi-invariant metric. Let $G^{\mathbb{C}}$ be the complexification
of $G$, with Lie algebra $\mathfrak{g}^{\mathbb{C}}$ (thus
$\mathfrak{g}^\mathbb{C}=\mathfrak{g}\otimes \mathbb{C}$).
Consider the based loop group
\begin{equation*}
\Omega G=\left\{\gamma:S^1\to G \,\,\mbox{(smooth)}\, |\,\,
\gamma(1)=e\right\}
\end{equation*}
and the corresponding infinite-dimensional Lie algebra
\begin{equation*}
\Omega \mathfrak{g}=\left\{\gamma:S^1\to \mathfrak{g}
\,\,\mbox{(smooth)}\, |\,\,
 \gamma(1)=0\right\}.
\end{equation*}

A smooth map $\Phi:\mathbb{C}\to \Omega G$ is called an
\emph{extended solution} if it satisfies
$$\Phi^{-1}d\Phi=\big(1-\lambda^{-1}\big)\alpha'+\big(1-\lambda\big)\alpha''$$
for each $\lambda\in S^1$, where $\alpha'$ is a $\g^\C$-valued
$(1,0)$-form on $\C$ with complex conjugate $\alpha''$. Observe
that, for each $z\in \C$, $\lambda\mapsto \Phi (z)(\lambda)$ is
holomorphic on $\C^*$.
\begin{thm}\label{uh}\cite{Uh}
\emph{\emph{a)} If $\Phi:\mathbb{C}\to \Omega G$ is an extended
solution, then the map $\phi:\mathbb{C}\to G$ defined by
$\phi(z)=\Phi (z)(-1)$ is harmonic. \emph{b)}  If
$\phi:\mathbb{C}\to G$ is harmonic, then there exists an extended
solution $\Phi:\mathbb{C}\to\Omega G$ such that $\phi(z)=\Phi
(z)(-1)$, for all $z\in \mathbb{C}$. This is unique up to
multiplication on the left by an element $\gamma\in\Omega G$ such
that $\gamma(-1)=e$.}
\end{thm}

\begin{rem}
Let $G/K$ be a symmetric space with automorphism $\tau$ and base
point $x_0=eK$. Define a map $\iota: G/K\to G$ by $\iota(g\cdot
x_0)=\tau(g)g^{-1}$. It is well known that $\iota$ is a totally
geodesic embedding, the \emph{Cartan embedding} of $G/K$ into $G$,
so that if $\varphi:\C\to G/K$ is harmonic then
$\iota\circ\varphi:\C\to G$ is also, and we can apply Theorem
\ref{uh} to $\iota\circ\varphi$.
\end{rem}

\section{Holomorphic potentials and extended framings}\label{si}

In \cite{DPW} the authors introduced a systematic procedure of
obtaining harmonic maps into a symmetric space from certain
holomorphic $1$-forms. The main ingredients are the existence of
various loop factorizations and the concept of \emph{extended
framing}. Next we will recall this construction.

\vspace{.15in}

Let $N=G/K$ be a symmetric space with automorphism $\tau$ and
associated symmetric decomposition $\g=\lk\oplus\mathfrak{m}$. Let
$\phi:\C\rightarrow N$ be a smooth map and take a lift
$\psi:\C\rightarrow G$ of $\phi$, that is, we have
$\phi=\pi\circ\psi$ where $\pi:G\rightarrow G/K$ is the coset
projection. Corresponding to the symmetric decomposition
$\g=\lk\oplus\mathfrak{m}$ there is a decomposition of
$\alpha=\psi^{-1}{d}\psi$, $
\alpha=\alpha_\mathfrak{k}+\alpha_\mathfrak{m}$.  Let
$\alpha_\mathfrak{m}=\alpha'_\mathfrak{m}+\alpha''_\mathfrak{m}$
be the type decomposition of $\alpha_\mathfrak{m}$ into
$(1,0)$-form and $(0,1)$-form of $\C$. Consider the loop of
$1$-forms
$\alpha_\lambda=\lambda^{-1}\alpha'_\mathfrak{m}+\alpha_\mathfrak{k}
+\lambda\alpha''_\mathfrak{m}$. We may view $\alpha_\lambda$ as a
$\Lambda_\tau\g$-valued $1$-form, where
\begin{equation}
\label{lambdatau}
 \Lambda_\tau\g=\left\{\xi:S^1\rightarrow
\g\,\,(\mathrm{smooth})\,|\,\,
\tau(\xi(\lambda))=\xi(-\lambda)\,\,\mathrm{for\,\,
all}\,\,\lambda\in S^1\right\}.
\end{equation}
It is well known that $\phi$ is harmonic if, and only if,
${d}+\alpha_\lambda$ is a loop of flat connections on the trivial
bundle $\underline{\C}^n=\C\times\C^n$. Hence, if $\phi$ is
harmonic, we can define a smooth map $\Psi:\C\rightarrow
\Lambda_\tau G$, where $\Lambda_\tau G$ is the
infinite-dimensional Lie group corresponding to the loop Lie
algebra (\ref{lambdatau}),
\begin{equation*}
\Lambda_\tau G=\left\{\gamma:S^1\rightarrow
G\,\,(\mathrm{smooth})\,|\,\,
\tau(\gamma(\lambda))=\gamma(-\lambda)\,\,\mathrm{for\,\,all}\,\,\lambda\in
S^1\right\},
\end{equation*}
such that $\Psi^{-1}d\Psi=\alpha_\lambda$. The smooth map $\Psi$
is called an \textit{extended framing} (associated to $\phi$). Our
harmonic map is recovered from $\Psi$ via $\phi=\pi\circ \Psi_1$
(here we are using the notation $\Psi_1(z)=\Psi (z)(1)$).

Consider the subspace $\Lambda_{-1,\infty}$ as in (\ref{prt}).
\begin{defn} \emph{A
holomorphic $1$-form $\mu$ on $\mathbb{C}$ with values in
$\Lambda_{-1,\infty}$ is called a \emph{holomorphic potential}.
The space of all holomorphic potentials is denoted by
$\mathcal{P}$. If $\mu\in\mathcal{P}$ takes values in
$\Lambda_{-1,\infty}^\tau =\Lambda_\tau\g^\C\cap
\Lambda_{-1,\infty}$, we say that $\mu$ is a \emph{$\tau$-twisted
holomorphic potential}. The space of all $\tau$-twisted
holomorphic potentials is denoted by $\mathcal{P}_\tau$. Thus
$\mu=\sum_{k\geqslant -1}\lambda^{k}\mu_k\in {\mathcal{P}}_\tau$
if $\mu_{even}$  is a $\lk^\C$-valued $1$-form and $\mu_{odd}$ is
a $\mathfrak{m}^\C$-valued $1$-form.}
\end{defn}

Fix an Iwasawa decomposition of the reductive group $K^\C$:
$K^\C=KB$ where $B$ is a solvable Lie subgroup. Consider the
following infinite-dimensional twisted Lie groups:
\begin{eqnarray*}
\Lambda_\tau G^\C&\!\!\!\!=\!\!\!\!&\left\{\gamma:S^1\rightarrow
G^\C\,\,(\mathrm{smooth})\,|\,
\tau(\gamma(\lambda))=\gamma(-\lambda)\,\,\mathrm{for\,\,all}\,\,\lambda\in
S^1\right\};\\
\Lambda^{+}_{B,\tau}
G^\C&\!\!\!\!=\!\!\!\!&\left\{\gamma\in\Lambda_\tau G^\C\,|\,
\gamma\,\, \mbox{extends holomorphically to
$|\lambda|<1$}\,,\,\gamma(0)\in B\right\}.
\end{eqnarray*}
We have the following twisted loop group decomposition:
\begin{thm}\label{tid}\cite{DPW}
\emph{Multiplication $\Lambda_\tau G\times \Lambda^+_{B,\tau}
G^\C\rightarrow \Lambda_\tau G^\C$ is a diffeomorphism onto.}
\end{thm}

So, let $\mu$ be a $\tau$-twisted holomorphic potential; since
$\mu$ is holomorphic, its $(0,1)$-part vanishes; then $\mu=\xi dz$
for some holomorphic function $\xi:\C\to \g^\C$ and it satisfies
the Maurer-Cartan equation $d\mu+\frac12[\mu\wedge\mu]=0$, that
is, $d_\mu=d+\mu$ is a flat connection. We can integrate to obtain
a unique map $\Psi_\mu:\C\to\Lambda_\tau G^\C$ such that
$\Psi_\mu^{-1}d\Psi_\mu=\mu$ and $\Psi_\mu(0)=e.$ If we factorize
$\Psi_\mu$ pointwise  according to Theorem \ref{tid}, we obtain a
smooth map $\Phi_\mu:\C\to \Lambda_\tau G$ such that
$\Psi_\mu=\Phi_\mu b,$ with $b:\C\to\Lambda^+_{B,\tau}G^C$. We
have:
\begin{thm}\cite{DPW} \emph{ $\Phi_\mu :\C\to \Lambda_\tau G$ is an extended
framing.}
\end{thm}
\section{Holomorphic potentials and extended solutions}\label{hol}

Since in our study of harmonic maps we will use extended solutions
instead of extended framings and, as in \cite{BP2}, we will
consider an action on the space of all extended solutions of germs
at zero of maps $\C\to G^\C$, in this section we reformulate the
DPW construction of harmonic maps according to our conveniences.
In particular, we have to enlarge the space of holomorphic
potentials. For completeness, we give detailed proofs of the main
results.

 \vspace{.15in}

Fix $0<\varepsilon <1$. Let $C_\varepsilon$ and
$C_{1/\varepsilon}$ denote the circles of radius $\varepsilon$ and
$1/\varepsilon$ centered at $0\in\mathbb{C}$; define open subsets
of $\mathbb{P}^1=\C\cup \{\infty\}$ by
$$I_\varepsilon=\left\{\lambda\in\mathbb{P}^1\,|\,\,|\lambda|<\varepsilon\right\},\,\,
I_{1/\varepsilon}=\left\{\lambda\in\mathbb{P}^1\,|\,\,|\lambda|>1/\varepsilon\right\},\,\,
E^\varepsilon=\left\{\lambda\in\mathbb{P}^1\,|\,\,\varepsilon<|\lambda|<1/\varepsilon\right\};$$
put $I^\varepsilon=I_\varepsilon\cup I_{1/\varepsilon}$ and
$C^\varepsilon=C_{\varepsilon}\cup C_{1/\varepsilon}$ so that
$\mathbb{P}^1=I^\varepsilon\cup C^\varepsilon\cup E^\varepsilon.$
Consider the infinite-dimensional Lie groups
\begin{align*}
\Lambda^\varepsilon G&=\big\{\gamma:C^\varepsilon\to
G^{\mathbb{C}}\,\, (\mathrm{smooth})\,|\,\,\overline{\gamma(\lambda)}=\gamma\left(1/\bar{\lambda}\right)\big\}\\
\Omega_E^\varepsilon G&=\left\{\gamma\in \Lambda^\varepsilon G
\,|\,\,\mbox{$\gamma$ extends holomorphically to $\gamma:
E^\varepsilon\to G^\C$ and $\gamma(1)=e$} \right\}\\
\Lambda_I^\varepsilon G&=\left\{\gamma\in \Lambda^\varepsilon G
\,|\,\,\mbox{$\gamma$ extends holomorphically to $\gamma:
I^\varepsilon\to G^\C$}\right\}
\end{align*}
and the corresponding infinite-dimensional Lie algebras
\begin{align*}
\Lambda^\varepsilon \mathfrak{g}&=\big\{\gamma:C^\varepsilon\to
\g^{\mathbb{C}}\,\, (\mathrm{smooth})\,|\,\,\overline{\gamma(\lambda)}=\gamma\left(1/\bar{\lambda}\right)\big\}\\
\Omega_E^\varepsilon \g&=\left\{\gamma\in \Lambda^\varepsilon \g
\,|\,\,\mbox{$\gamma$ extends holomorphically to $\gamma:
E^\varepsilon\to \g^\C$ and  $\gamma(1)=0$} \right\}\\
\Lambda_I^\varepsilon \g &=\big\{\gamma\in \Lambda^\varepsilon \g
\,|\,\,\mbox{$\gamma$ extends holomorphically to $\gamma:
I^\varepsilon\to \g^\C$}\big\}.
\end{align*}

To express the extended solution equation in terms of
$\Lambda^\varepsilon G$ we shall make  use of the following
Iwasawa-type decomposition of $\Lambda^\varepsilon G$:
\begin{thm}\cite{Mc}\label{Mc}
\emph{Multiplication $\Omega^\varepsilon_E G \times
\Lambda^\varepsilon_IG \to\Lambda^\varepsilon G$ is a
diffeomorphism. In particular, any $\gamma\in \Lambda^\varepsilon
G$ may by written uniquely in the form $\gamma=\gamma_E\gamma_I$,
where $\gamma_E\in \Omega^\varepsilon_E G$ and $\gamma_I\in
\Lambda^\varepsilon_I G^\mathbb{C}.$}
\end{thm}

\begin{rem}\label{Moca}
The limiting case of Theorem \ref{Mc} as $\varepsilon\rightarrow
1$ is the more familiar decomposition $\Omega G\times
\Lambda_+G^\mathbb{C}\to\Lambda G^\mathbb{C}$ where $\Lambda
G^{\mathbb{C}}=\left\{\gamma:S^1\to
G^{\mathbb{C}}\,|\,\,\mbox{$\gamma$ is smooth}\right\}$ and
$$\Lambda_+ G^{\mathbb{C}}=\left\{\gamma\in \Lambda
G^\C\,|\,\,\mbox{$\gamma$ extends holomorphically to $
|\lambda|<1$ }\right\}.$$ This result is due to Pressely-Segal
\cite{PS}.
\end{rem}

We also have to enlarge our class of potentials: fix
$0<\varepsilon < 1$ and consider the subspace of
$\Lambda^\varepsilon\g$ defined by
$\Lambda^\varepsilon_{-1,\infty}=\left\{\xi\in
\Lambda^\varepsilon\g\,|\,\,\mbox{$\lambda \xi$ extends
holomorphically to $I_\varepsilon$}\right\};$ each element $\xi\in
\Lambda^\varepsilon_{-1,\infty}$ can be written as
$\xi=(\xi_+,\xi_{-})$, where $\xi_+:C_\varepsilon \to \g^\C$
extends meromorphically to $I_\varepsilon$ with at most a simple
pole at $0$ and $\xi_{-}:C_{1/\varepsilon}\to \g^\C$ is defined by
$\xi_{-}(\lambda)=\overline{\xi_+\big(1/\bar{\lambda}\big)};$
\begin{defn}\label{pot} \emph{
 A $1$-form $\mu=(\mu_+,\mu_-)$ on $\mathbb{C}$ with values in
 $\Lambda^\varepsilon_{-1,\infty}$ such that $\mu_+$ is
 holomorphic is called a \emph{$\varepsilon$-holomorphic potential}. The space of
all $\varepsilon$-holomorphic potentials is denoted by
$\mathcal{P}^\varepsilon$.}
\end{defn}
\begin{rem}
 The space of holomorphic
potentials $\mathcal{P}$ can be interpreted as the limiting case
of Definition \ref{pot} as $\varepsilon\rightarrow 1$.
\end{rem}


 Let $\mu$ be a holomorphic
$\varepsilon$-potential, so that $d_\mu=d+\mu$ is a flat
connection. This means that we can integrate to obtain a unique
map $\Psi_\mu:\C\to\Lambda^\varepsilon G$, with
$\Psi_{\mu}^{-1}d\Psi_\mu = \mu$  and  $\Psi_\mu (0)=e.$ We call
$\Psi_\mu$ a \emph{complex extended solution}. If we factorize
$\Psi_\mu$ according to Theorem \ref{Mc}, we obtain a map
$\Phi_\mu:\C\to \Omega^\varepsilon_E G$ such that
$\Psi_\mu=\Phi_\mu b$,
 with $b:\C\to \Lambda^\varepsilon_I G$. Note that, since $\Psi_\mu(0)=e$,
 we also have $\Phi_\mu(0)=e$ and $b(0)=e$.
\begin{thm}\label{no} \emph{ $\Phi_\mu:\C\to\Omega^\varepsilon_E G\subset \Omega G $  is a extended
solution.}
\end{thm}
\begin{proof}With respect to the Iwasawa decomposition, the Lie algebra
$\Lambda^\varepsilon \mathfrak{g}$ splits into a direct sum of Lie
subalgebras:
\begin{equation}\label{split}
\Lambda^\varepsilon\mathfrak{g}=\Omega^\varepsilon_E
\mathfrak{g}\oplus \Lambda^\varepsilon_I \mathfrak{g}.
\end{equation}
Since $\Phi_\mu=\Psi_\mu b^{-1}$, we have
\begin{equation}\label{derivada}
\Phi_\mu^{-1}d \Phi_\mu= \mathrm{Ad}_b(\mu)-db\, b^{-1}.
\end{equation}
But $b$ takes values in $\Lambda^\varepsilon_I G$, that is
\begin{equation}\label{b}
b(z)=b_0(z)+b_1(z)\lambda+b_2(z)\lambda^{2}+\ldots
\end{equation}
for all $\lambda\in I_\varepsilon$,
 so that $db\, b^{-1}$ takes values in $\Lambda_I^\varepsilon
\g$, then by (\ref{split}) and (\ref{derivada})
\begin{equation}
\Phi_\mu^{-1}d
\Phi_\mu=\big(\mathrm{Ad}_b(\mu)\big)_{\Omega^\varepsilon_E \g}.
\end{equation}
Now, $\mu$  is  a $1$-form on $\mathbb{C}$ with values in
$\Lambda^\varepsilon_{-1,\infty}$ and $\Lambda_I^\varepsilon G$
acts on $\Lambda^\varepsilon_{-1,\infty}$; in $C_\varepsilon$ we
can write $\mu=\sum_{k\geqslant -1}\mu_k\lambda^{k}$;  hence
\begin{equation}
\big(\mathrm{Ad}_b(\mu)\big)_{\Omega^\varepsilon_E
\g}=\big(\lambda^{-1}-1\big)\mathrm{Ad}_{b_0}(\mu_{-1})+\big(\lambda-1\big)\overline{\mathrm{Ad}_{b_0}(\mu_{-1})}.
\end{equation}
Since $\mu_{-1}$ is a $\g^\C$-valued holomorphic $1$-form  we can
write $\mu_{-1}= \xi_{-1}dz$ for some holomorphic function
$\xi_{-1}:\C\to \g^\C$; hence
\begin{equation}\label{der}
\Phi_\mu^{-1}d
\Phi_\mu=\big(1-\lambda^{-1}\big)\alpha'+\big(1-\lambda\big)\alpha'',
\end{equation}
with
\begin{equation}\label{ad}
\alpha'=-\mathrm{Ad}_{b_0}(\mu_{-1})=-\mathrm{Ad}_{b_0}(\xi_{-1})dz,\end{equation}
that is, $\Phi_\mu$ is an extended solution.
\end{proof}
Then any $\varepsilon$-holomorphic potential $\mu$ gives rise to a
harmonic map $\phi_\mu: \C \to G$ with
$\phi_\mu(z)=\Phi_\mu(z)(-1)$
 and $\phi_\mu(0)=\Phi_\mu(0)(-1)=e$.
 \begin{rem}
Again, by taking the limiting case of Theorem \ref{no} as
$\varepsilon\to 1$ we see that any holomorphic potential $\mu \in
\mathcal{P}$ gives rise to a harmonic map $\phi_\mu:\C\to G$.
 \end{rem}

\begin{rem}\label{rem}
Set $\alpha_\lambda=(1-\lambda^{-1})\alpha'+(1-\lambda)\alpha''$.
Let $d=\partial+\bar{\partial}$ and
$d_{\alpha_\lambda}=\partial_{\alpha_\lambda}+\bar{\partial}_{\alpha_\lambda}$,
respectively, be the type decompositions of the connections $d$
(the trivial connection) and
$d_{\alpha_\lambda}=d+\alpha_\lambda$, respectively, on
$\underline{\C}^{n}=\C\times\C^n$. Considering the $(0,1)$-parts
of equations (\ref{derivada}) and (\ref{der}), we obtain
$-(1-\lambda)\alpha''=\bar{\partial}b\,b^{-1}.$ This means that
the gauge transformation $b$ gauges the holomorphic structure
$\bar{\partial}_{\alpha_\lambda}$ on $\underline{\C}^{n}$  to the
trivial holomorphic structure $\bar{\partial}$ on
$\underline{\C}^{n}$, that is,
$b^{-1}\bar{\partial}_{\alpha_\lambda}b=\bar{\partial}.$ In
particular, $b:(\underline{\C}^{n},\bar{\partial})\to
(\underline{\C}^{n},\bar{\partial}_{\alpha_\lambda})$ is a
holomorphic isomorphism.
\end{rem}

\begin{rem}\label{pita} Denote by $\hat{\mathcal{P}}$ the space of all
$1$-forms on $\C$ satisfying the following conditions: \emph{a)}
$\mu$ takes values in $\Lambda_{-1,\infty}$; \emph{b)}
$d_\mu=d+\mu$ is a flat connection; \emph{c)} the $(0,1)$-part of
$\mu$ takes values in $\Lambda_+\g^{\C}$, the Lie algebra of
$\Lambda_+G^\C$. Again we can integrate and apply the Iwasawa
decomposition of Remark \ref{Moca}
 to obtain a smooth map $\Phi_\mu:\C\to\Omega G$. Moreover, the
arguments given to prove Theorem \ref{no} carry over directly to
this case to prove that $\Phi_\mu$ is also an extended solution.
This  space of potentials will be useful for us in section
\ref{sopas}.
\end{rem}

Any harmonic map from $\C$ to $G$ is obtained, up to left
multiplication by a constant, from a holomorphic potential
$\mu\in\mathcal{P}$:
\begin{thm}
\emph{Let $\phi:\C\to G$ be a harmonic map such that $\phi(0)=e$.
Then there exists a holomorphic potential $\mu\in\mathcal{P}$ such
that $\phi=\phi_\mu$.} \end{thm}
\begin{proof}
Let $\phi:\C\to G$ be a harmonic map such that $\phi(0)=e$. Let
$\Phi$ be an extended solution associated to $\phi$, with
$\Phi^{-1}d\Phi=\big(1-\lambda^{-1}\big)\alpha'+\big(1-\lambda\big)\alpha''$.
Consider  the $\Lambda_+\g^\C$-valued $(0,1)$-form
$\theta=-(1-\lambda)\alpha''$. The $\bar{\partial}$-problem
$\bar{\partial}b\,b^{-1}=\theta,b(0)=e $ can be solved over $\C$
(see appendix of \cite{DPW}). Let $b:\C\to\Lambda_+G^{C}$ be the
unique solution of this problem and put $\Psi=\Phi b$. We have
$\Psi^{-1}d\Psi=\mathrm{Ad}_{b^{-1}}\big(\Phi^{-1}d\Phi\big)+b^{-1}db.$
By construction we see that $\mu=\Psi^{-1}d\Psi$ is a
$\Lambda_{-1,\infty}$-valued $1$-form of type $(1,0)$ (in
particular, due to Maurer-Cartan equation, $\mu$ is holomorphic)
such that $\phi=\phi_\mu$.
\end{proof}

\section{Extended framings vs. extended solutions}\label{bb} In the
previous sections we have seen two procedures of obtaining
harmonic maps from a potential $\mu$: via extended solutions we
get an harmonic map $\phi:\C\to G$; via extended framings (if
$\mu$ is a twisted potential) we get an harmonic map
$\tilde{\phi}:\C\to G/K$. We show in this section that (when $\mu$
is a twisted potential) $\phi$ and $\tilde{\phi}$ are essentially
the same harmonic map.

\vspace{.15in}

Fix a twisted potential $\mu\in{\mathcal{P}}_\tau$ and integrate
to obtain a complex extended solution $\Psi:\C\to\Lambda G^\C$.
This map has a unique factorization $\Psi=\Phi b$, with
$\Phi:\C\to \Omega G$ an extended solution and $b:\C\to \Lambda_+
G^\C$. The corresponding harmonic map into $G$ is given by
$\phi=\Phi_{-1}$ (again, we are using the notation
$\Phi_{\lambda}(z)=\Phi(z)(\lambda)$). On the other hand, we can
view $\Psi$ as a map from $\C$ to $\Lambda_\tau G^\C$ and use the
decomposition of Theorem \ref{tid} to write
$\Psi=\tilde{\Phi}\tilde{b}$, where
$\tilde{\Phi}:\C\to\Lambda_\tau G$ is an extended framing and
$\tilde{b}:\C\to \Lambda^+_{B,\tau}G^\C$. The corresponding
harmonic map into the symmetric space   $G/K$ is given by
$\tilde{\phi}=\pi\circ\tilde{\Phi}_1$. Let $\iota:G/K\to G$ be the
Cartan embedding. It happens that $\phi$ has values in $\iota
(G/K)$ and $\phi=\iota\circ\tilde{\phi}$. In fact: we have
$\iota\circ\pi\circ\tilde{\Phi}_{1}=\tau(\tilde{\Phi}_1)\tilde{\Phi}_1^{-1}=\tilde{\Phi}_{-1}\tilde{\Phi}_1^{-1};$
however, by the uniqueness of the decomposition $\Psi=\Phi b$, we
have $\Phi=\tilde{\Phi}\tilde{\Phi}_1^{-1}$ and
$b=\tilde{\Phi}_1\tilde{b}$, so that
$\iota\circ\pi\circ\tilde{\Phi}_{1}=\tilde{\Phi}_{-1}\tilde{\Phi}_1^{-1}=\Phi_{-1};$
hence $\Phi$ and $\tilde{\Phi}$ produce the same harmonic map into
 $G/K$.

\section{Dressing actions and gauge transformations}

Another consequence of Iwasawa-type decomposition of Theorem
\ref{Mc} is that it allows us to define a natural action
$\#_\varepsilon$ of $\Lambda^\varepsilon_I G$ on
$\Omega^\varepsilon_E G$: if $g\in\Omega^\varepsilon_E G$ and
$h\in \Lambda^\varepsilon_I G$, then $h\#_\varepsilon g=(hg)_E.$
Applying this action pointwise, we obtain from an extended
solution $\Phi:\C\to \Omega G$ a new map $h\#_\varepsilon\Phi:
\C\to \Omega G$ defined by
$(h\#_\varepsilon\Phi)(z)=h\#_\varepsilon\Phi(z)$, for all
$z\in\C$.

\begin{thm}\cite{BP2,DPW,GO}
\emph{$h\#_\varepsilon\Phi$ is an extended solution.}
\end{thm}
We will recall now from \cite{BP2} how these actions vary with
$\varepsilon$. For $0<\varepsilon<\varepsilon'<1$ we have
injections $\Lambda_I^{\varepsilon'}G\subset
\Lambda_I^{\varepsilon}G$ and $\Omega_E^{\varepsilon}G\subset
\Omega_E^{\varepsilon'}G$. Similarly, for $0<\varepsilon<1$, we
have $\Omega_{\rm{hol}}G\subset \Omega_E^{\varepsilon}G$, where
$\Omega_{\rm{hol}}G=\bigcap_{0<\varepsilon <1}
\Omega_E^\varepsilon G$. Its is easy to see that
\begin{align*}
\Omega_{\rm{hol}}G & =\big\{\gamma:\C^*\to
G^\C\,|\,\,\mbox{$\gamma$ is holomorphic, $\gamma(1)=e$ and
  $\overline{\gamma(\lambda)}=\gamma\big(1/\bar{\lambda}\big)$}
\big\}.
\end{align*}
We have:
\begin{thm}\cite{BP2}
\emph{For $0<\varepsilon<\varepsilon'<1$,
$\gamma\in\Lambda_I^{\varepsilon'}G \subset
\Lambda_I^{\varepsilon}G$, and $g\in\Omega_E^\varepsilon
G\subset\Omega_E^{\varepsilon'}G$, we have
$\gamma\#_{\varepsilon'} g=\gamma\#_\varepsilon g\in
\Omega^\varepsilon_E G.$}
\end{thm}
\begin{axiom}\cite{BP2}
\emph{The action of each $\Lambda_I^{\varepsilon'}G$ preserves
$\Omega_{\rm{hol}}G$ and, for $0<\varepsilon<\varepsilon'<1,$
$\gamma\in\Lambda_I^{\varepsilon'}G \subset
\Lambda_I^{\varepsilon}G$ and $g\in\Omega_{\rm{hol}}G$, we have
$\gamma\#_{\varepsilon'} g=\gamma\#_\varepsilon g$.}
\end{axiom}
It follows that we can take a direct limit as $\varepsilon\to 0$
and so obtain an action on $\Omega_{\rm{hol}}G$ of the group of
germs at zero of maps $\C\to G^\C$. Henceforth, we write $\gamma\#
g$ for this action on $\Omega_{\rm{hol}}G$.

On the other hand, the holomorphic gauge group
$$\mathcal{G}^\varepsilon=\big\{\mbox{$h=(h_+,h_-):\C\to
\Lambda^\varepsilon_IG\,\,$ such that
$\bar{\partial}h_+=0$}\big\}$$ acts on the space
$\mathcal{P}^\varepsilon$ of $\varepsilon$-holomorphic potentials
by gauge transformations: if $\mu\in \mathcal{P}^\varepsilon$ and
$h\in\mathcal{G}^\varepsilon$, then
$h\cdot\mu=\mathrm{Ad}_h(\mu)-dh\,h^{-1}\in
\mathcal{P}^\varepsilon$. It happens that the correspondence
between holomorphic potentials and extended solutions
$\mu\to\Phi_\mu$ is equivariant with respect to these actions:
\begin{thm}\label{roi}\cite{BP2,DPW}
\emph{If $h\in\mathcal{G}^\varepsilon$, then
$\Phi_{h\cdot\mu}=h(0)\#\Phi_\mu.$}
\end{thm}
\begin{rem}
The limiting case of Theorem \ref{roi} as $\varepsilon\to 1$ can
be stated as follows:

Consider the dressing action $\#$ of $\Lambda_+G^\C$ on $\Omega G$
corresponding to the Iwasawa decomposition of Remark \ref{Moca};
let $\mathcal{G}$ be the holomorphic gauge group of all $h:\C\to
\Lambda _+ G^\C$ such that $\bar{\partial}h=0$; let $\mu\in
\mathcal{P}$ be a holomorphic potential. Then
$\Phi_{h\cdot\mu}=h(0)\#\Phi_\mu$.
\end{rem}
\section{Adding a uniton and gauge
transformations}\label{sopas}

Another well known operation for generating new extended solution
from a given one was introduced by Uhlenbeck \cite{Uh} and is
called {\it{adding a uniton}}. In this section we describe how
this operation arises via gauge transformations on the holomorphic
potential. More precisely: consider a holomorphic potential
$\mu\in\mathcal{P}$ with associated extended solution $\Phi_\mu$;
denote by $S_\mu$ the set of all gauge transformations $h:\C\to
\Lambda G^\C$ such that $h\cdot \mu=\mathrm{Ad}_h(\mu)-dh\,h^{-1}$
is in $\hat{\mathcal{P}}$ (see Remark \ref{pita}); each element
$h$ of $S_\mu$ gives rise to a new extended solution
$\Phi_{h\cdot\mu}$; in particular, we have $\mathcal{G}\subset
S_\mu$ and $\Phi_{h\cdot\mu}=h(0)\#\Phi_\mu$
 if $h\in \mathcal{G}$; next we describe the elements of $h\in
S_\mu$ whose action on $\mu$ corresponds to the operation of
adding a uniton to $\Phi_\mu$.

\vspace{.15in}

\begin{thm}\cite{Uh}\label{adduni}
\emph{Let $\Phi:\C\rightarrow \Omega {U}(n)$ be an extended
solution, $\phi:\C\rightarrow {U}(n)$ the corresponding harmonic
map. Write $\alpha_0=\frac12
\phi^{-1}d\phi=A_zdz+A_{\bar{z}}d\bar{z}$. Let $\hat{\ell}$ be a
subbundle of $\,\underline{\C}^{n}$ with Hermitian projection
$\hat{\pi}: \underline{\C}^{n}\rightarrow \hat{\ell}$ satisfying
the \emph{{uniton conditions}}: \emph{a)} $\hat{\pi}^\perp
A_z\hat{\pi} =0$; \emph{b)}
$\hat{\pi}^\perp(\bar{\partial}\hat{\pi}+
A_{\bar{z}}\hat{\pi})=0.$
 Then
$\tilde{ \Phi}:\C\rightarrow \Omega {U}(n)$ given by $\tilde{
\Phi}_\lambda=\Phi_\lambda(\hat{\pi}+\lambda^{}\hat{\pi}^\perp)$
is an extended solution.} \end{thm}

This operation of obtaining new harmonic maps from a given one is
called \emph{adding a uniton} in \cite{Uh}. Note that the second
uniton condition means that $\hat{\ell}$ is a holomorphic
subbundle of $\underline{\C}^{n}$ with respect to the holomorphic
structure $\bar{\partial}_{\alpha_0}$.

\vspace{.05in}

 We  recall from \cite{Uh} how to add a
uniton to a harmonic map into a Grassmannian:

Let   $G_k(\C^n)$ be the complex Grassmannian of $k$-planes in
$\C^n$. The unitary group   ${U}(n)$ acts transitively on
$G_k(\C^n)$ with stabilizers conjugate to ${ U}(k)\times{
U}(n-k)$. Fix a complex $k$-plane  $V_0\in G_k(\C^n)$ with
stabilizer $K$ and let  $\pi_0$ be the Hermitian projection onto
$V_0$. Let $\tau$ be the involution of ${U}(n)$ given by
conjugation by $Q_0=\pi_0-{\pi_0}^\perp$. The identity component
of the fixed set of $\tau$ is $K$ so that $G_k(\C^n)$ is a
symmetric space with involution $\tau$. The corresponding Cartan
embedding $\iota_k: G_k(\C^n) \rightarrow {U}(n)$ is given by $
\iota_k(V)=Q_0(\pi_V-{\pi_V}^\perp$), where $\pi_V$ denotes the
Hermitian projection onto the $k$-plane $V$.

\begin{thm}
\cite{Uh}\label{addunigr} \emph{Suppose that $\psi:\C\rightarrow
G_k(\C^n)$ is a harmonic map and
 $\Phi$ an extended solution associated to $\phi=\iota_k\circ \psi$.
 Let $\hat{\ell}$ be a subbundle of
  $\underline{\C}^{n}$, with Hermitian projection
 $\hat{\pi}: \underline{\C}^{n}\rightarrow \hat{\ell}$, satisfying the {{uniton
 conditions}}, and  such that
$[\phi,\hat{\pi}]=0$.
 Then $\tilde{\Phi}:\C\rightarrow \Omega {U}(n)$ given by
$\tilde{
\Phi}_\lambda=\Phi_\lambda(\hat{\pi}+\lambda\hat{\pi}^\perp)$ is
an extended solution associated to a harmonic map $\tilde{\psi}$
into a Grassmannian  $G_{\tilde{k}}(\C^n)$:
$\tilde{\Phi}_{-1}=\iota_{\tilde{k}}\circ \tilde{\psi}$.}
\end{thm}

Let $\ell$ be a holomorphic subbundle of
$(\underline{\C}^{n},\bar{\partial})$, with Hermitian projection
$\pi$. Fix a holomorphic potential $\mu\in\mathcal{P}$. Suppose
that $\pi^\perp\mu_{-1}\pi=0$. Define $\gamma_\ell:\C\to\Omega G$
by
\begin{equation}\label{gammaell}
\gamma_\ell=\pi+\lambda^{-1} \pi^\perp.
\end{equation}
This map $\gamma_\ell$ gauges the connection $d_\mu$ to the
connection associated to the $1$-form
$$\gamma_\ell\cdot \mu =\left(\pi+\lambda^{-1}
\pi^\perp\right)\Big(\sum_{k\geqslant
-1}\mu_k\lambda^k\Big)\left(\pi+\lambda
\pi^\perp\right)-\left(1-\lambda^{-1}\right)d\pi \left(\pi+\lambda
\pi^\perp\right).$$ The coefficient in $\lambda^{-2}$ on the right
hand of this equality is zero since $\pi^\perp\mu_{-1}\pi=0$.
Moreover, since $\ell$ is a holomorphic subbundle of
$\underline{\C}^n$, we have $\bar{\partial} \pi\pi=0$; hence the
$(0,1)$-part of the coefficient in $\lambda^{-1}$ on the right
hand of this equality also vanishes. Then
$\gamma_\ell\cdot\mu\in\hat{\mathcal{P}}$ and $\gamma_\ell\in
S_\mu$.
\begin{thm}\label{pt}
\emph{$\Phi_{\gamma_\ell\cdot \mu}$ is obtained from $\Phi_\mu$ by
adding a uniton. More precisely: consider the Iwasawa
decomposition of $\Psi_\mu$, $\Psi_\mu=\Phi_\mu b$, and the
subbundle $\hat{\ell}=b_0\ell$, where $b_0(z)=b(z)(0)$, with
Hermitian projection $\hat{\pi}: \underline{\C}^{n}\rightarrow
\hat{\ell}$; then $\hat{\ell}$ satisfies the uniton conditions of
Theorem \ref{adduni} and $$\Phi_{\gamma_\ell\cdot
\mu}=(\hat{\pi}_0+\lambda^{-1}\hat{\pi}_0^\perp)\Phi_\mu(\hat{\pi}+\lambda\hat{\pi}^\perp),$$
where $\hat{\pi}_0$ is the Hermitian projection onto the fibre of
$\hat{\ell}$ at $z=0$.}
\end{thm}
\begin{proof}
First note that
$\Psi_{\gamma_\ell\cdot\mu}=(\hat{\pi}_0+\lambda^{-1}\hat{\pi}_0^\perp)\Psi_\mu({\pi}+\lambda{\pi}^\perp)$,
i.e.,
$\Psi^{-1}_{\gamma_\ell\cdot\mu}d\Psi_{\gamma_\ell\cdot\mu}=\gamma_\ell\cdot\mu$
and $\Psi_{\gamma_\ell\cdot\mu}(0)=e.$ Then, to find
$\Phi_{\gamma_\ell\cdot\mu}$ we have to factorize
$\Psi_{\gamma_\ell\cdot\mu}$ according to the Iwasawa
decomposition. We have
\begin{align}\label{coco}
\Psi_{\gamma_\ell\cdot\mu}&=(\hat{\pi}_0+\lambda^{-1}\hat{\pi}_0^\perp)\Psi_\mu({\pi}+\lambda{\pi}^\perp)=
(\hat{\pi}_0+\lambda^{-1}\hat{\pi}_0^\perp)\Phi_\mu
b({\pi}+\lambda{\pi}^\perp)\nonumber
\\&= (\hat{\pi}_0+\lambda^{-1}\hat{\pi}_0^\perp)\Phi_\mu
(\hat{\pi}+\lambda\hat{\pi}^\perp)(\hat{\pi}+\lambda^{-1}\hat{\pi}^\perp)
b({\pi}+\lambda{\pi}^\perp).
\end{align}
We claim that $\hat{b}=(\hat{\pi}+\lambda^{-1}\hat{\pi}^\perp)
b({\pi}+\lambda{\pi}^\perp)$ takes values in $\Lambda_+ G^\C$. In
fact, $\hat{b}$ is holomorphic at $\lambda=0$ if $\hat{\pi}^\perp
b_0\pi=0$ (this is the coefficient in $\lambda^{-1}$). But
$\hat{\pi}^\perp b_0\pi=0$  if and only if $\hat{\ell}=b_0\ell$,
whence  $\hat{\pi}^\perp b_0\pi$ vanishes automatically. The
invertibility follows by applying the same argument to
$\hat{b}^{-1}$. Whence $\hat{b}$ takes values in $\Lambda_+ G^\C$.
From (\ref{coco}) we  conclude now that $\Phi_{\gamma_\ell\cdot
\mu}=(\hat{\pi}_0+\lambda^{-1}\hat{\pi}_0^\perp)\Phi_\mu(\hat{\pi}+\lambda\hat{\pi}^\perp)$
(in particular, $\hat{\ell}$ satisfies the uniton conditions.)
\end{proof}
Reciprocally, any uniton can be added via the action of some loop
of the form (\ref{gammaell}) on the holomorphic potential $\mu$:
\begin{thm}
\emph{Suppose that the subbundle $\hat{\ell}$ with Hermitian
projection
 $\hat{\pi}: \underline{\C}^{n}\rightarrow \hat{\ell}$ satisfies the uniton conditions
 with respect to $\Phi_\mu$. Set $\ell=b_0^{-1}\hat{\ell}$. Then:
\emph{a)} $\ell$ is holomorphic with respect to the trivial
holomorphic structure $\bar{\partial}$; \emph{b)}
$\pi^\perp\mu_{-1}\pi=0$ and
 \emph{c)} we have $\Phi_{\gamma_\ell\cdot
\mu}=(\hat{\pi}_0+\lambda^{-1}\hat{\pi}_0^\perp)\Phi_\mu(\hat{\pi}+\lambda\hat{\pi}^\perp).$}
\end{thm}
\begin{proof}
\emph{a)} We have seen in Remark \ref{rem} that
$b:(\underline{\C}^n,\bar{\partial})\to
(\underline{\C}^n,\bar{\partial}_{\alpha_\lambda})$ is a
holomorphic isomorphism; hence, $\ell$ is a holomorphic subbundle
of $(\underline{\C}^n,\bar{\partial})$, since $\hat{\ell}$ is a
holomorphic subbundle of
$(\underline{\C}^n,\bar{\partial}_{\alpha_0})$. \emph{b)} Equation
$\pi^\perp\mu_{-1}\pi=0$ follows directly from the first uniton
condition and (\ref{ad}) . \emph{c)} This statement follows
directly from Theorem \ref{pt}.
\end{proof}
\section{Harmonic maps into Grassmannians and subbundles}
 As an application of the ideas of section
\ref{sopas},  in section \ref{s} we will be able to improve and
clarify  some results presented in the unpublished second author's
doctoral thesis \cite{Pa} concerning unitons preserving
\emph{finite type} property of harmonic maps. Before we shall
recall from \cite{BW} some relevant facts about harmonic maps into
Grassmannians:

\vspace{.15in}

Let $N=G/K$ be a symmetric space with involution $\tau$. Denote by
$\g$ and $\lk$ the Lie algebras of $G$ and $K$, respectively, and
consider the corresponding symmetric decomposition
$\g=\lk\oplus\m$ into $\pm 1$-eigenspaces of the derivation of
$\tau$. Recall that, for each $x=g\cdot x_0$, the surjective map
$\g\to T_x N$ given by $\xi\mapsto \frac{d}{dt}{\big|_{t=0}} \exp
t\xi\cdot x$ has the Lie algebra ${\rm{Ad}}_g\lk$ as kernel and so
restricts to an isomorphism ${\rm{Ad}}_g \m \to T_xN$. The inverse
map $\beta_x:T_xN\to {\rm{Ad}}_g\m$ defines a $\g$-valued $1$-form
$\beta$ on $N$, the \emph{Maurer-Cartan form} of $N=G/K$. We
denote by $[\m]$ the subbundle of the trivial bundle
$\underline{\g}=N\times \g$ defined by $[\m]_{g\cdot
x_0}=\mathrm{Ad}_g(\m)$.

 If $N$ is actually the (matrix) group manifold $G$, acting on
itself by right translations, then $\beta$ is just the (left)
Maurer-Cartan form $\theta$ of $G$. Moreover:
\begin{lem}\cite{BR}\label{pullback} \emph{If $\psi:\C\to G/K$ is a smooth map
and $\iota$ is the Cartan embedding of $G/K$ into $G$, then, with
$\phi=\iota\circ\psi$ we have
$(\phi^{-1}d\phi=)\phi^*\theta=-2\psi^*\beta.$}
\end{lem}
\begin{proof}
In fact: let $X\in T_xG/K$; then $X=\frac{d}{dt}\big|_{t=0}\exp
t\beta(X)\cdot x$ so that
\begin{align*}
\iota_*(X)&=\frac{d}{dt}\Big|_{t=0}\tau\big(\exp
t\beta(X)\big)\iota(x)\exp
\big(-t\beta(X)\big)\\&=\iota(x)\big(\mathrm{Ad}_{\iota(x)^{-1}}\tau(\beta(X))-\beta(X)\big)=-2\iota(x)\beta(X);
\end{align*}
hence $\iota^*\theta=\iota^{-1}d\iota=-2\beta$, and from this
equality, by taking the pullback by $\psi$, we see that the
statement of this lemma holds.
\end{proof}

Consider now the case of the Grassmannian $G_k(\C^n)$. Denote by
$\mathfrak{u}(n)$ the Lie algebra of $U(n)$ and take $V_0\in
G_k(\C^n)$ as base point. The associated symmetric decomposition
$\mathfrak{u}(n)=\mathfrak{k}\oplus\mathfrak{m}$ is given by
\begin{eqnarray*}
\mathfrak{k}^\C&=&{\rm{Hom}}(V_0,V_0)\oplus{\rm{Hom}}({V_0}^\perp,{V_0}^\perp)\,,
\\\mathfrak{m}^\C&=&{\rm{Hom}}(V_0,{V_0}^\perp)\oplus{\rm{Hom}}({V_0}^\perp,V_0)\,.
\end{eqnarray*}
Let $T\rightarrow G_k(\C^n)$ be the tautological subbundle of
$G_k(\C^n)\times \C^n$ whose fibre at $V\in G_k(\C^n)$  is $V$
itself. With respect to the usual Hermitian structure of
$G_k(\C^n)$, the Maurer-Cartan form of $G_k(\C^n)$ restricted to
the bundle of $(1,0)$-vectors gives an isomorphism
$\beta^{(1,0)}:T^{(1,0)}G_k(\C^n)\rightarrow
{\rm{Hom}}(T,T^\perp)$.

 Identify  a smooth map $\psi:\C\rightarrow G_k(\C^n)$ with
the smooth complex subbundle $\underline{\psi}$ of the trivial
bundle $\underline{\C}^{n}$ given by setting the fibre at $z$
equal to $\psi(z)$ for all $z\in \C$. Conversely any rank $k$
subbundle of  $\underline{\C}^{n}$ induces a map $\C\rightarrow
G_k(\C^n)$. Denote the Hermitian  projection onto a vector
subbundle $\underline{\psi}$ by $\pi_{\psi}$ and define vector
bundle morphisms $A'_{\psi},A''_{\psi}:
\underline{\psi}\rightarrow\underline{\psi}^\perp$ called the
$\partial$- and $\bar{\partial}$-\emph{second fundamental forms}
of $\underline{\psi}$ in  $\underline{\C}^n$ by
$A'_{\psi}(v)=\pi_{\psi^\perp}(\frac{\partial v}{\partial z})$ and
$A''_{\psi}(v)=\pi_{\psi^\perp}(\frac{\partial v}{\partial
\bar{z}})$, for $v$  a smooth section of $\underline{\psi}$. Note
that $A'_\psi$ is minus the adjoint of $A''_{\psi^\perp}$. The
second fundamental forms of $\underline{\psi}$ in
$\underline{\C}^n$, $A'_\psi$ and $A''_\psi$, represent, via
$\beta^{(1,0)}$, the $(1,0)$-components of the partial derivatives
$\partial \psi$ and $\bar{\partial}\psi$ (see \cite{BW}):
\begin{equation}\label{pb}
\psi^*\beta\Big(\frac{\partial}{\partial
z}\Big)=\psi^*\beta^{(1,0)}\Big(\frac{\partial}{\partial
z}\Big)+\psi^*\beta^{(0,1)}\Big(\frac{\partial}{\partial
z}\Big)=A'_\psi+A'_{\psi^\perp}.
\end{equation}
Comparing Lemma \ref{pullback} with (\ref{pb}) we obtain the
following formula to the derivative of $\phi=\iota\circ\psi$ in
terms of the second fundamental forms of $\psi$:
\begin{equation}\label{ves}
\frac12\phi^{-1}\partial\phi=-(A'_\psi+A'_{\psi^\perp}).
\end{equation}

In this setting, the harmonicity equations can be reformulated as
follows (cf. \cite{BW}): we give each subbundle of
$\underline{\C}^{n}$ the connection induced from the trivial
connection of $\underline{\C}^{n}$ and corresponding
Koszul--Malgrange holomorphic structure; a smooth map
$\psi:\C\rightarrow G_k(\C^n)$ is harmonic if and only if
$A'_\psi$ is holomorphic; this holds if and only if $A''_\psi$  is
anti-holomorphic.

Now, given two holomorphic bundles $E$, $F$ over a Riemann surface
and  a holomorphic bundle morphism $A:E \rightarrow F$,  there are
unique holomorphic subbundles $\underline{\mathrm{Ker}}A$ and
$\underline{\mathrm{Im}}A$ of $E$ and $F$, respectively, that
coincide with  ${\mathrm{Ker}}A$ and ${\mathrm{Im}}A$ almost
everywhere. Thus we can  define the $\partial$- and
$\bar{\partial}$-\emph{Gauss bundles} of a harmonic map
$\psi:\C\to G_k(\C^n)$  by
$G^{(1)}(\psi)=\underline{\mathrm{Im}}A'_{\psi}$ and
$G^{(-1)}(\psi)=\underline{\mathrm{Im}}A''_{\psi}$, respectively.
These bundles also represent harmonic maps (cf. \cite{BW}).
Iterating this construction we set
$G^{(0)}(\psi)=\underline{\psi}$, and for $i=1,2,\ldots$,
$G^{(i)}(\psi)=G'(G^{(i-1)}(\psi)),
G^{(-i)}(\psi)=G''(G^{-(i-1)}(\psi)).$ The bundle $G^{(i)}(\psi)$
is called {\it{the $i^{th}$-Gauss bundle of $\psi$}}.

\begin{rem}\label{rema}
The harmonic map $G^{(-1)}(\psi)$ can be obtained by adding the
uniton $\ell=\underline{\ker} A'_{\psi^\perp}$ to $\psi$.
\end{rem}

\section{Unitons preserving  finite type}\label{s}

\begin{defn}\cite{BP}\emph{ A harmonic map $\phi:\C\to G$ is of \emph{finite type} if it
can be obtained from an extended solution $\Phi_\mu:\C\to\Omega G$
whose holomorphic potential $\mu=\xi dz$ is constant of the form
$\xi=\lambda^{d-1}\eta$, for some odd $d\in\mathbb{N}$, with
$$\eta\in\Omega_d\g=\big\{\eta\in\Omega\g\,|\,\,\eta=\sum_{|k|\leq
d}\eta_k\lambda^{k}\big\}.$$}
\end{defn}
It is clear that:
\begin{thm}\label{ha}
\emph{Fix a harmonic map $\phi_\mu:\C\to G$ of finite type, with
$\mu=\lambda^{d-1}\eta dz$ and $\eta\in\Omega_d\g$. Fix a constant
subspace of $\C^n$ and let $\pi_0:\C^n\to \ell_0$ be the
corresponding hermitian projection. If
$\gamma_{\ell_0}=\pi_0+\lambda^{-1}\pi_0^\perp$ is a  loop in
$S_\mu$, then $\phi_{\gamma_{l_0}\cdot\mu}$ is also of finite
type. }\end{thm} Consider the vector subspace $\ell_0=\ker
\eta_{-d}$ of $\C^n$, with Hermitian projection $\pi_0$; in this
case, we have $\pi_0^\perp\mu_{-1}\pi_0=0$, whence
$\gamma_{\ell_0}\in S_\mu$.
 A more concrete example of unitons preserving finite type is given by the following
theorem:

\begin{thm}\emph{ If $\psi:\C\to G_k(\C^n)$ is of finite type, then, for
each integer $r$, $G^{(r)}(\psi)$ is of finite type.}
\end{thm}
\begin{proof}Fix $\psi(0)$ as base point of $G_k(\C^n)$,
let $K$ be the isotropic subgroup of $U(n)$ associated to
$\psi(0)$ and $\tau$ the corresponding automorphism. Let $\mu=\xi
dz\in\mathcal{P}_\tau$ be a  constant $\tau$-twisted holomorphic
potential of the form $\xi=\lambda^{d-1}\eta$, for some odd
$d\in\mathbb{N}$, with $\eta=\sum_{|k|\leq d}\eta_k\lambda^{k}\in
\Omega_d \g,$ associated to $\psi$.
 Since $d$ is odd
 and $\eta$ is twisted,
we have $\eta_{-d}\in \mathfrak{m}^\C$.  According to the
decomposition of $\m^\C$ into $(1,0)$-vectors and $(0,1)$-vectors,
$\mathfrak{m}^\C=\mathfrak{m}^+\oplus\mathfrak{m}^-,$ with
$\mathfrak{m}^+=\mathrm{Hom}(\psi(0),\psi(0)^\perp)$ and
$\mathfrak{m}^-=\mathrm{Hom}(\psi(0)^\perp,\psi(0)),$ write
$\eta_{-d}=\eta^+_{-d}+\eta^-_{-d}$. Fix the vector subspace
$\ell_0=\ker \eta^-_{-d}$. Clearly we have $\gamma_{\ell_0}\in
S_\mu$; thus $\Phi_{\gamma_{\ell_0}\cdot \mu}$ gives rise to a
harmonic map of finite type. We claim that this harmonic map is
precisely the Gauss bundle $G^{-1}(\psi)$:

Let $\Psi_\mu$ be the complex extended solution associated to
$\mu$. Factorize $\Psi_\mu$ according to Remark \ref{Moca},
$\Psi_\mu=\Phi b$.  Fix an Iwasawa decomposition $K^\C=KB$ and
factorize $\Psi_\mu$ according to Theorem \ref{tid},
$\Psi_\mu=\tilde{\Phi} \tilde{b}$. Since \emph{i)} $\tilde{b}_0$
takes values in $B\subset K^\C$, \emph{ii)} $\tilde{\Phi}_1$
frames $\psi$, that is, $\psi=\tilde{\Phi}_1\cdot \psi(0)$, and
\emph{iii)} $b_0=\tilde{\Phi}_1\tilde{b}_0$ (cf. Section
\ref{bb}), we conclude that $b_0:\C\to Gl(n,\C)$ also frames
$\psi$, that is, $\psi=b_0\cdot \psi(0)$. In particular,
$\psi^*[\mathfrak{m}^\C]_z=\mathrm{Ad}_{b_0(z)}(\mathfrak{m}^\C)$
and
\begin{equation}\label{hom}
\mathrm{Hom}(\psi^*T,\psi^*T^\perp)_z=\mathrm{Ad}_{b_0(z)}(\mathfrak{m}^+)\quad\quad
\mathrm{Hom}(\psi^*T^\perp,\psi^*T)_z=\mathrm{Ad}_{b_0(z)}(\mathfrak{m}^-)
\end{equation}
for each $z\in\C$. On the other hand, from (\ref{ad}) and
(\ref{ves}) follows that
\begin{equation}
\label{porra}
 A'_\psi+A'_{\psi^\perp}=\mathrm{Ad}_{b_0}(\eta_{-d}).
\end{equation}
Hence, from (\ref{hom}) and (\ref{porra}) we conclude that
$A'_{\psi^\perp}=\mathrm{Ad}_{b_0}(\eta^-_{-d})$. Thus,
$$\underline{\ker}\,
A'_{\psi^\perp}=\underline{\ker}\,\mathrm{Ad}_{b_0}(\eta^-_{-d})={b_0}\ell_0.$$
Taking account Remark \ref{rema} and Theorem \ref{pt}, this
establishes the claim.

Iterating this argument and reversing orientation we have the
result.
\end{proof}

However, unitons do not always preserve finite type: let
$\phi:\C\to G_k(\C^n)$ be a (non-constant) harmonic map of finite
type and $\delta:\C\to G_s(\C^m)$ a (non-constant) holomorphic
map. Then $\psi=\phi\oplus \delta :\C\to G_{k+s}(\C^n\oplus\C^m)$
is a harmonic map which is obtained from $\phi$ by adding the
uniton $\delta$. If $A'_\delta$ has singular points, the same
happens to $A'_\psi=A'_\phi\oplus A'_\delta$; in this case, $\psi$
can not be of finite type, since equation (\ref{porra}) ensures
that the second fundamental forms $A'_\psi$ and $A''_{\psi^\perp}$
associated to a harmonic map $\psi:\C\to G_k(\C^n)$ of finite type
have no singular points (points where the rank of $\mathrm{Im}
A'_\psi$ drops).

\begin{rem}Let $\{\gamma_a\}$ be a curve in $\Lambda^\varepsilon_IG$ and
$\Phi$  an extended solution. Suppose $\lim_{a\to
0}\gamma_a=\gamma$ with $\gamma \in\Lambda^\varepsilon G$. The
extended solution $\tilde{\Phi}=\lim_{a\to 0}\big(\gamma_a\#
\Phi\big)$ need not to coincide with $(\gamma \Phi)_E$; thus,
$\tilde{\Phi}$ is not, apriori, obtained by applying a dressing
transformation to $\Phi$. Following \cite{bG}, we shall refer to
the procedure of obtaining $\tilde{\Phi}$ from $\Phi$ as
\emph{completion}. In \cite{Uh}, Uhlenbeck suggested that any
uniton can be added this way. Bergvelt and Guest \cite{bG} settled
negatively this conjecture and studied in detail this procedure of
completion in the case of curves in $\Lambda^\varepsilon_IG$
defined by $\gamma_a=\pi_V+\xi_a\pi_V^\perp,$
 where $V$ is a constant subspace of $\C^n$, $a\in\C$, and
$\xi_a$ is the rational function in $\lambda$ given by
\begin{equation}\label{sp}
\xi_a(\lambda)=\frac{\bar{a}\lambda-1}{\lambda-a}\frac{1-a}{\bar{a}-1}.
\end{equation}
We observe that if $\Phi_\mu$ is an extended solution associated
to the holomorphic potential $\mu\in\mathcal{P}$, and $V$ is a
constant subspace of $\C^n$ such that $\pi_V^\perp\mu_{-1}\pi_V=0$
everywhere, then $\lim_{a\to 0}\gamma_a\cdot \mu =\gamma\cdot
\mu\in{\mathcal{P}}$, with $\gamma=\pi_V+\lambda^{-1}\pi_V^\perp$.
Hence, in this case, the completion procedure amounts to add a
uniton. In particular, we see that the unitons preserving finite
type in Theorem \ref{ha} can be added via completion.

\end{rem}

%


\begin{thebibliography}{10}
\bibitem{bG} M.J. Bergvelt and M.A. Guest, \emph{Action of loop groups
on harmonic maps.} Trans. Amer. Math. Soc. \textbf{326} (1991),
861--886
\bibitem{BFP} F.E. Burstall, D. Ferus, F. Pedit, and U.
Pinkall, \emph{{H}armonic Tori in symmetric spaces and commuting
Hamiltonian systems on loop algebras}, Ann. of Math. \textbf{138}
(1993), 173--212.
\bibitem{BP2}
F.E. Burstall and F. Pedit, \emph{{D}ressing orbits of harmonic
maps}, Duke Math. J. \textbf{80} (1995), 353--382.
\bibitem{BP}
F.E. Burstall and F. Pedit, \emph{{H}armonic maps via
Adler-Konstant-Symes theory}, Harmonic maps and Integrable Systems
(A.P. Fordy and J.C.Wood, eds), Aspects of Mathematics E23,
Vieweg, 1994, pp. 221--272. CMP 94:09
\bibitem{BR}
F.E. Burstall, J. H. Rawnsley, \emph{{T}wistor Theory for
Riemannian Symmetric Spaces}, Lectures Notes in Math. 1424 Berlin,
Heidelberg: 1990.
\bibitem{BW}
F.E. Burstall and J.C. Wood, \emph{{T}he construction of harmonic
maps into complex Grassmannians}, J. Diff. Geom. \textbf{23}
(1986), 255--297.

\bibitem{DPW}
J. Dorfmeister, F. Pedit and H. Wu, \emph{Weiestrass type
representation of harmonic maps into symmetric spaces}, Comm.
Anal. Geom. \textbf{6} (1998), 633-668.

%
\bibitem{GO}
M.A. Guest and Y. Ohnita, \emph{{G}roup actions and deformations
for harmonic maps}, J. Math. Soc. Japan \textbf{45} (1993),
671--710.
\bibitem{Mc} I. McIntosh, \emph{Global solutions of the elliptic 2D
 periodic Toda lattice,} Nonlinearity \textbf{7} (1994), 85-108.

\bibitem{Pa}
R. Pacheco, \emph{Harmonic maps and loop groups}, Ph.D. thesis,
University of Bath, 2004.



\bibitem{PS}
A.N. Presseley and G.B. Segal, \emph{{L}oop Groups}, Oxford
University Press, 1986.
\bibitem{Uh}
K. Uhlenbeck, \emph{{H}armonic maps into Lie groups (classical
solutions of the chiral model)}, J. Diff. Geom. \textbf{30}
(1989), 1--50.
\end{thebibliography}
\end{document}